\documentclass[review]{elsarticle}

\usepackage{bm,amssymb,amsmath,mathrsfs}
\usepackage{amsthm}
\usepackage{lineno}
\usepackage[usenames]{color}

\usepackage{dcolumn}

\newcommand{\bq}{\begin{equation}}
\newcommand{\eq}{\end{equation}}
\newcommand{\bc}{\begin{center}}
\newcommand{\ec}{\end{center}}
\newcommand{\bit}{\begin{itemize}}
\newcommand{\eit}{\end{itemize}}
\newcommand{\ben}{\begin{enumerate}}
\newcommand{\een}{\end{enumerate}}

\theoremstyle{plain}

\newtheorem*{theorem*}{Theorem}

\begin{document}

\journal{(internal report CC24-1)}

\begin{frontmatter}

\title{Probability Mass Function, Moments and Factorial Moments of the Negative Binomial Distribution NB$(k,r)$}

\author[cc]{S.~R.~Mane}
\ead{srmane001@gmail.com}
\address[cc]{Convergent Computing Inc., P.~O.~Box 561, Shoreham, NY 11786, USA}

\begin{abstract}
  The negative binomial distribution NB$(k,r)$ of Type I is the probability distribution for a sequence of independent Bernoulli trials
  (with success parameter $p\in(0,1)$) with $r$ nonoverlapping success runs of length $\ge k$.
  We present a new, more concise, expression for its probability mass function.
  We show it can also be succinctly written using hypergeometric functions.
  We also present new expressions (combinatorial sums) for its moments and factorial moments, as opposed to only the mean and variance (which are already known).
  Next, we present an alternative non-combinatorial viewpoint, which yields expressions for the factorial moments not only for nonoverlapping success runs,
  but also for runs with an overlap of $\ell$, where $\ell\in[0,k-1]$.
  The case $\ell=k-1$ is the negative binomial distribution NB$(k,r)$ of Type III.
  The results also yield the solution for the negative binomial distribution NB$(k,r)$ with a minimum gap between the success runs (explained in the text).
  $\mathit{Addendum\; 1/23/2024:}$ The probability mass function and factorial moments are derived from the probability generating function.
  $\mathit{Addendum\; 1/26/2024:}$ Alternative expressions are presented for the negative binomial distribution NB$(k,r)$ of Type II.
\end{abstract}

\vskip 0.25in


\end{frontmatter}

\newpage
\setcounter{equation}{0}
\section{\label{sec:intro}Introduction}
The negative binomial distribution NB$(k,r)$ (or NB${}_{k,r}$) of Type I is the probability distribution for a sequence of $n$ independent Bernoulli trials
(all with success parameter $p\in(0,1)$)
with $r$ nonoverlapping success runs of length $\ge k$.
The term ``nonoverlapping'' follows the definition in the classic text by Feller \cite{Feller}.
This means that the counting for a new success run begins as soon as a success run of length $k$ is attained.
Hence a run of six consecutive successes $SSSSSS$ signifies
(i) six consecutive success runs if $k=1$,
(ii) three consecutive success runs if $k=2$,
(iii) two consecutive success runs if $k=3$,
(iv) one success run if $k=6$
and for $k=4$ and $k=5$ it contains one success run, with the beginnings of another success run.
This is called the negative binomial distribution of Type I.
The text by Balakrishnan and Koutras \cite{BalakrishnanKoutras} gives an excellent exposition on the subject,
with formulas for the negative binomial distribution of Types I, II and III.
\begin{enumerate}
\item
The negative binomial distribution of Type II requires at least one failure between success runs.
In the above example, the sequence $SSSSSS$ counts as only one success run, for any value $k=1,\dots,6$.
We do not treat Type II in this note.
\item
The negative binomial distribution of Type III permits $k-1$ overlaps between consecutive success runs of length $k$.
In the above example, the sequence $SSSSSS$ contains
(i) two success runs of length 5 if $k=5$,
(ii) three success runs of length 4 if $k=4$,
(iii) four success runs of length 3 if $k=3$,
(iv) five success runs of length 2 if $k=2$.
We discuss Type III briefly in Sec.~\ref{sec:DM4}.
\end{enumerate}
Note that \cite{BalakrishnanKoutras} is a reference text:
our principal references for the original formulas are 
Charalambides \cite{Charalambides1986} and Philippou and Georghiou \cite{PhilippouGeorghiou1989}.
Both sets of authors published expressions for the probability mass function as nested combinatorial sums.
Charalambides \cite{Charalambides1986} also published a recurrence for the factorial moments.

This note presents a new, more concise, expression for the probability mass function (pmf) of the negative binomial distribution NB${}_{k,r}$ of Type I.
The pmf can also be succinctly written using hypergeometric functions (which is computationally faster to evaluate).
We also present new expressions (combinatorial sums) for its moments and factorial moments, as opposed to only the mean and variance. 
Next, we present an alternative non-combinatorial viewpoint to derive the factorial moments.
It applies also for runs with an overlap of $\ell$, where $\ell\in[0,k-1]$.
Finally, the results also yield the solution for the negative binomial distribution NB$(k,r)$ with a minimum gap between the success runs (see Sec.~\ref{sec:gap}).

The structure of this paper is as follows.
Sec.~\ref{sec:notation} presents the basic definitions and notation employed in this note.
Secs.~\ref{sec:pmf}, \ref{sec:facmom} and \ref{sec:mom} present our results for the probability mass function, factorial moments and moments, respectively.
Sec.~\ref{sec:cent} gives a brief discussion of the central moments,
i.e.~moments centered on the mean, in particular the skewness and kurtosis.
Sec.~\ref{sec:DM4} presents an alternative non-combinatorial expression for the factorial moments.
Sec.~\ref{sec:gap} presents the application of the above results to the case of the negative binomial distribution NB$(k,r)$ with a minimum gap between the success runs.
Sec.~\ref{sec:conc} concludes.

\newpage
\setcounter{equation}{0}
\section{\label{sec:notation}Basic notation and definitions}
For future reference, we denote the probability mass function (pmf)
of the negative binomial distribution NB$(k,r)$ of Type I by $P_{k,r}(n)$.
In most of this note we shall drop the dependence on $k$ and $r$ and write $P_n$ for brevity.
(This is also the notation employed in \cite{Charalambides1986} and \cite{PhilippouGeorghiou1989}.)
We also denote its moments and factorial moments by $M_n$ and $M_{(n)}$, respectively.
{\em Do not be confused by the different uses of $n$:}
in $P_n$ it denotes the number of Bernoulli trials,
but $M_n = \mathbb{E}[Y^n]$,
and $M_{(n)} = \mathbb{E}[(Y)_n] = \mathbb{E}[Y(Y-1)\dots(Y-n+1)]$,
where $Y \sim \mathrm{NB}_{k,r}$.
The use of $n$ in this way is because we run out of letters of the alphabet otherwise:
see the notation in the formulas in the published literature.

For future reference, define $q=1-p$.
First consider the {\em geometric distribution of order $k$}.
This is the distribution of the waiting time for the {\em first} run of $k$ consecutive successes, hence its pmf is $P_{k,1}(n)$, i.e.~$r=1$.
It has been extensively studied \cite{Feller}--\cite{DM3}.
As pointed out in \cite{BalakrishnanKoutras} and \cite{PhilippouGeorghiou1989},
the negative binomial distribution NB$(k,r)$ of Type I is obtained via an $r$-fold convolution of identical geometric distributions of order $k$.
As mentioned above, we denote its pmf by $P_n$.
Charalambides \cite{Charalambides1986} and Philippou and Georghiou \cite{PhilippouGeorghiou1989}
published (equivalent) recurrence relations for $P_n$.
Both sets of authors also published expressions for $P_n$ as nested combinatorial sums.
Charalambides \cite{Charalambides1986} also published a recurrence relation for the factorial moment for $M_{(n)}$, but not an explicit solution.

To avoid confusion, we clarify the indexing of trials.
Many authors, such as Feller \cite{Feller}, begin by counting the first Bernoulli trial as $n=1$.
Philippou and Georghiou \cite{PhilippouGeorghiou1989} employ the same notation.
In that case, we require at least $rk$ trials to obtain $r$ nonoverlapping success runs of length $\ge k$, hence $P_n=0$ for $n \in [1,rk-1]$ and $P_n = p^{rk}$ for $n=rk$.
However, other authors, such as Charalambides \cite{Charalambides1986}, begin the indexing when it is possible to attain a positive probability and define $P_0 = p^{rk}$ and begin the indexing of $n$ from there.
To avoid confusion, we employ the terms `$F$' and `$C$' indexing (for `full' and `cut') to distinguish the two definitions of indexing.
{\em Note that this is not standard terminology.}
With an obvious notation, the relation is $n_C = n_F -rk$.

\newpage
With this understanding of the indexing, the recurrence for the pmf by Philippou and Georghiou is
(Theorem 3.1 in \cite{PhilippouGeorghiou1989}), where $n$ refers to $F$ indexing
\bq
\label{eq:Pn_PG_Thm3.1}
P_n = \begin{cases} 0 &\qquad n < rk \,,
  \\
  p^{rk} &\qquad n=rk \,,
  \\
  \displaystyle \frac{q/p}{n-rk}\sum_{j=1}^k (n-rk +j(r-1))p^jP_{n-j} &\qquad n > rk \,.
  \end{cases}
\eq
An equivalent recurrence was published by Charalambides (eq.~(4.17) in \cite{Charalambides1986})\footnote{Some misprints in \cite{Charalambides1986} have been corrected.}
(where $n$ to refers to $C$ indexing)
\bq  
\label{eq:Pn_Ch_eq4.17}
P_n = \frac{q/p}{n}\sum_{j=1}^{\min(n,k)} (n +rj-j)p^jP_{n-j} \,.
\eq
Charalambides also published a recurrence for the factorial moment $M_{(n)}$ (eq.~(4.18) in \cite{Charalambides1986}), where $n$ to refers to $C$ indexing
\bq
\label{eq:Mn_rec_Ch}
\begin{split}
  M_{(n)} &= \theta(1-\theta)^{-1} \sum_{j=1}^{\min(n,k)} \binom{n-1}{j-1}\frac{n+rj-j}{j}\mu_{(j)}M_{(n-j)}
  \\
  &= \frac{1-p^k}{np^k} \sum_{j=1}^{\min(n,k)} \binom{n}{j}(n+rj-j)\mu_{(j)}M_{(n-j)} \,.
\end{split}
\eq
Here $\theta = 1-p^k$ and 
$\mu_{(j)} = \mathbb{E}[(X)_j]$ (unnumbered, after eq.~(4.3) in \cite{Charalambides1986}),
where $X\in[1,k]$ and $P(X=j) = (1-p^k)^{-1}qp^{j-1}$
(unnumbered, before eq.~(2.1) in \cite{Charalambides1986}).
See Sec.~\ref{sec:facmom} for further processing of $\mu_{(j)}$.

We employ eqs.~\eqref{eq:Pn_Ch_eq4.17} and \eqref{eq:Mn_rec_Ch} as our principal formulas below (with $C$ indexing),
then transcribe/transform the results to $F$ indexing.

Before proceeding further, we comment on a result by Muselli,
who treated the negative binomial distribution of Type II with $r$ success runs of length $\ge k$ and obtained the following expression for the pmf,
where the formula used $x$ instead of $n$ and employed $F$ indexing (Theorem 4 in \cite{Muselli}).
\bq
\label{eq:Muselli_pmf}
P(x) = \sum_{m=r}^{\lfloor\frac{x+1}{k+1}\rfloor} (-1)^{m-r}\binom{m-1}{r-1}p^{mk}q^{m-1}\biggl(\binom{x-mk-1}{m-2} +q\binom{x-mk-1}{m-1}\biggr) \,.
\eq
Although we shall not treat the negative binomial distribution of Type II in this note, Muselli's formula is important.
It gives us a clue as to the form of the answer we seek.
We wish to obtain a similarly elegant expression (to the extent possible) for the negative binomial distribution of Type I,
for the pmf and then the moments and factorial moments.
That is the goal of this note.

\newpage
\setcounter{equation}{0}
\section{\label{sec:pmf}Probability mass function}
\subsection{$C$ indexing}
We begin with $C$ indexing.
Charalambides published an explicit sum for the pmf (eq.~(4.19) in \cite{Charalambides1986})\footnote{Some misprints in \cite{Charalambides1986} have been corrected.}
\bq
\label{eq:Ch_sum_full}  
\begin{split}
P_n &= \sum_{i=0}^n \binom{r+i-1}{r-1} \,q^ip^{n+rk-i} \biggl[ \sum_{j=0}^i (-1)^j\binom{i}{j}\binom{n-jk-1}{i-1} \biggr] \,.
\end{split}
\eq
Note the following.
The sum over $i$ terminates at $n-jk$ else the last binomial coefficient is zero.
The sum over $j$ terminates when $n-jk-1 < 0$.
Hence set $n-jk-1=0$ to deduce $j \le \lfloor(n-1)/k\rfloor$.
Using eq.~\eqref{eq:Muselli_pmf} as a guide to organize the calculation,
reorganize eq.~\eqref{eq:Ch_sum_full} to nest the sum with $j$ on the outside and $i$ on the inside and simplify.
The result is 
\bq
\label{eq:pmf_sum_C}
\begin{split}
  P_n = p^{rk} \sum_{j=0}^{\lfloor(n-1)/k\rfloor} (-1)^j p^{jk} q^j\binom{r+j-1}{r-1} \biggl[\sum_{i=0}^r \binom{r}{i} \binom{n-jk-1}{j+i-1} q^i\biggr] \,.
\end{split}
\eq
Compared to eq.~\eqref{eq:Ch_sum_full}, there are fewer terms to sum.
For fixed $r$ and $k$, eq.~\eqref{eq:Ch_sum_full} contains $(n+1)(n+2)/2 = O(n^2)$ terms to sum.
In eq.~\eqref{eq:pmf_sum_C}, the value of $j$ goes up to $\lfloor(n-1)/k\rfloor$.
For fixed $j$, the number of terms to sum is not $2$ as in eq.~\eqref{eq:Muselli_pmf},
but $r+1$ (which is however independent of $n$).
Hence for fixed $r$ and $k$, an upper bound on the number of terms to sum in eq.~\eqref{eq:pmf_sum_C} is $(r+1)(1+\lfloor(n-1)/k\rfloor) = O(n)$.
\begin{enumerate}
\item
  Actually, eq.~\eqref{eq:pmf_sum_C} contains binomial coefficients which equal zero,
so there are fewer terms to sum.
\item
  When $j=0$, we must exclude the first term ($i=0$) because $\binom{n-1}{-1}=0$.
\item
  We could set the lower limit to $i=\delta_{j,0}$.
\item
  At the top, the sum over $i$ cuts off when $j+i-1 = n-jk-1$, i.e.~$i = n-(j+1)k$.
\item
  Instead of fussing with complicated upper and lower limits, it is more elegant to leave
eq.~\eqref{eq:pmf_sum_C} `as is' and skip terms where the binomial coefficient is zero.
\end{enumerate}
The sum over $i$ in eq.~\eqref{eq:pmf_sum_C} can be expressed more concisely as a hypergeometric function.
The answer is
\bq
\label{eq:pmf_hyp_C}
\begin{split}
  P_n &= p^{rk} \biggl[ qr\, {}_2F_1(1-n,1-r,2;q)
    \\
    &\qquad\qquad
    +\sum_{j=1}^{\lfloor(n-1)/k\rfloor} (-1)^j p^{jk} q^j\binom{r+j-1}{r-1}\binom{n-jk-1}{j-1} {}_2F_1(jk+j-n,-r,j;q)
    \biggr] \,.
\end{split}
\eq
The term in $j=0$ must be treated as a special case because we cannot set $j=0$ in the denominator of a hypergeometric function.
This is because the first binomial coefficient ($i=0$) is zero when $j=0$ in eq.~\eqref{eq:pmf_sum_C}.
The expression in eq.~\eqref{eq:pmf_hyp_C} is less intuitive to interpret or vizsualize, but is more efficient to compute.

\newpage
\subsection{$F$ indexing}
For $F$ indexing, we replace $n \gets n-rk$. Then
\bq
\label{eq:pmf_sum_F}
\begin{split}
  P_n = p^{rk} \sum_{j=0}^{\lfloor(n-rk-1)/k\rfloor} (-1)^j p^{jk} q^j\binom{r+j-1}{r-1} \biggl[\sum_{i=0}^r \binom{r}{i} \binom{n-(r+j)k-1}{j+i-1} q^i\biggr] \,.
\end{split}
\eq
Using hypergeometric functions,
\bq
\label{eq:pmf_hyp_F}
\begin{split}
  P_n &= p^{rk} \biggl[ qr\, {}_2F_1(rk+1-n,1-r,2;q)
    \\
    &\qquad\qquad
    +\sum_{j=1}^{\lfloor(n-rk-1)/k\rfloor} (-1)^j p^{jk} q^j\binom{n-(r+j)k-1}{j-1}\binom{r+j-1}{r-1} {}_2F_1((r+j)k+j-n,-r,j;q)
    \biggr] \,.
\end{split}
\eq

\newpage
\setcounter{equation}{0}
\section{\label{sec:facmom}Factorial moments}
\subsection{General remarks}
We first note that the calculation of the moments and factorial moments has a similarity
to the corresponding calculation for the Poisson distribution of order $k$.
Charalambides published a recurrence for the factorial moments of 
the Poisson distribution of order $k$ (eq.~(4.7) in \cite{Charalambides1986}),
which has a broadly similar structure to eq.~\eqref{eq:Mn_rec_Ch}.
Charalambides also published a solution (combinatorial sum) for the factorial moments of 
the Poisson distribution of order $k$ (eq.~(4.9) in \cite{Charalambides1986}).
That solution was independently obtained by the author \cite{Mane_CC23_11}
and in a subsequent note \cite{Mane_CC23_12} the methodology was extended to solve for the moments of the Poisson distribution of order $k$.
The solution below for the moments and factorial moments of the negative binomial distribution NB$(k,r)$ of Type I
follows the ideas in \cite{Mane_CC23_11,Mane_CC23_12}.

\subsection{$C$ indexing}
We begin with $C$ indexing.
The equation to solve is eq.~\eqref{eq:Mn_rec_Ch} and we require an expression for $\mu_{(j)}$.
It is the $j^{th}$ factorial moment of the random variable $X$ in the text below eq.~\eqref{eq:Mn_rec_Ch}.
\bq
\begin{split}
\mu_{(j)} &= \frac{q}{1-p^k} \sum_{i=1}^k(i)_jp^{i-1}
\\
&= \frac{q}{1-p^k} \sum_{i=j}^k i(i-1)\dots(i-j+1)p^{i-1} 
\\
&\equiv \frac{C_jp^k}{1-p^k} \,.
\end{split}
\eq
Some algebra yields
\bq
\label{eq:Cj_facmom}
C_j = \frac{j!}{q^jp^k}\biggl[ p^{j-1}(1-p^k) -p^k\sum_{i=1}^{j-1} \binom{k}{i} p^{j-i-1}q^i \biggr] -(k)_j \,.
\eq
Here $(k)_j = k(k-1)\dots(k-j+1)$.
Substitute into eq.~\eqref{eq:Mn_rec_Ch} to obtain
\bq
  M_{(n)} = \frac{1}{n} \sum_{j=1}^{\min(n,k)} \binom{n}{j}(n+rj-j)C_jM_{(n-j)} \,.
\eq
Following the ideas in \cite{Mane_CC23_11,Mane_CC23_12}, the solution is as follows.
Define $s_n = n_1+\dots+n_k$. Then
\bq  
\label{eq:facmom_C}
M_{(n)} = n! \sum_{n_1+2n_2+\dots+kn_k=n} (r+s_n-1)_{s_n} \prod_{j=1}^k \frac{1}{n_j!}\Bigl(\frac{C_j}{j!}\Bigr)^{n_j} \,.
\eq

\newpage
\subsection{$F$ indexing}
We cannot simply shift the factorial means by (powers of) $rk$.
Let $Y \sim \mathrm{NB}_{k,r}$ with $F$ indexing and $Z \sim \mathrm{NB}_{k,r}$ with $C$ indexing.
Then, writing $a=rk$ for brevity,
\bq
\begin{split}
  Y_{(j)} &= (Z+a)_{(j)}
  \\
  &= (Z+a)(Z-1+a)\dots(Z-j+1+a)
  \\
  &= Z_{(j)} +\binom{j}{1}Z_{(j-1)}a_{(1)} +\binom{j}{2}Z_{(j-2)}a_{(2)} +\dots +a_{(j)}
  \\
  &= \sum_{i=0}^j \binom{j}{i} Z_{(j-i)}a_{(i)} \,.
\end{split}
\eq
Denote the factorial moment for $F$ indexing by $\hat{M}_{(n)}$.
It follows that
\bq
\label{es:facmom_shift}
\hat{M}_{(n)} = \sum_{i=0}^j \binom{j}{i} M_{(n-i)}a_{(i)} \,.
\eq
It remains to determine an explicit expression for $\hat{M}_{(n)}$.
Analogously to eq.~\eqref{eq:Cj_facmom}, define $F_j$ as follows.
Some algebra with $C_j$ and shifts yields
\bq
\label{eq:Fj_facmom}
F_j = \frac{j!}{q^jp^k}\biggl[ p^{j-1}(1-p^k) +\sum_{i=1}^{j-1} (-1)^i \binom{k+i-1}{i} p^{j-i-1}q^i \biggr]
\eq
Unlike the expression for $C_j$ in eq.~\eqref{eq:Cj_facmom}, there is no tail piece $(k)_j$.
Also the terms in the sum alternate in sign and are not multiplied by $p^k$
and the binomial coefficient in the sum is different.
Then, using $F$ indexing for $n$, we obtain the following result.
Recall $s_n = n_1+\dots+n_k$.
\bq  
\label{eq:facmom_F}
\hat{M}_{(n)} = n! \sum_{n_1+2n_2+\dots+kn_k=n} (r+s_n-1)_{s_n} \prod_{j=1}^k \frac{1}{n_j!}\Bigl(\frac{F_j}{j!}\Bigr)^{n_j} \,.
\eq
Both eqs.~\eqref{eq:facmom_C} and \eqref{eq:facmom_F} have the same structure, just substitute $C_j$ for $C$ indexing and $F_j$ for $F$ indexing.

\newpage
\subsection{Mean}
For $n=1$ we obtain the obvious shift of the mean by $rk$.
\bq
\hat{M}_{(1)} = M_{(1)} + a_{(1)} = M_{(1)} + rk \,.
\eq
It is well known that the expressions for the mean, say $\mu_C$ and $\mu_F$ for $C$ and $F$ indexing respectively, are
\begin{subequations}
\begin{align}
\mu_C &= rC_1 \quad = r\,\frac{1-p^k}{qp^k} -rk \,,
\\
\mu_F &= rF_1 \quad = r\,\frac{1-p^k}{qp^k} \,.
\end{align}
\end{subequations}
It is evident that $\mu_F = \mu_C +rk$.

\newpage
\subsection{Variance}
For $n=2$ we obtain the following relation.
\bq
\begin{split}
\hat{M}_{(2)} &= M_{(2)} +2M_{(1)}a_{(1)} +a_{(2)}
\\
&= M_{(2)} +2rkM_{(1)} +rk(rk-1) \,.
\end{split}
\eq
The variance $\sigma^2$ is the same using either indexing, because it does not depend on the origin of the counting of the Bernoulli trials.
\begin{enumerate}  
\item
  Using $C$ indexing yields
\bq
\label{eq:var_facmom_C}  
\sigma^2 = M_{(2)} - M^2_{(1)} +M_{(1)} \,.
\eq
\item
  Using $F$ indexing yields
\bq
\label{eq:var_facmom_F}  
\begin{split}
  \sigma^2 &= \hat{M}_{(2)} - \hat{M}^2_{(1)} +\hat{M}_{(1)}
  \\
  &= M_{(2)} +2rkM_{(1)} +rk(rk-1) -(M_{(1)} + rk)^2  +M_{(1)} + rk
  \\
  &= M_{(2)} -M_{(1)}^2 +M_{(1)} \,.
\end{split}
\eq
\item
  Both eqs.~\eqref{eq:var_facmom_C} and \eqref{eq:var_facmom_F} yield the same expression.
\item
  The variance is proportional to $r$ and has no term in $r^2$. We can confirm this as follows.
\begin{subequations}
\begin{align}
M_{(1)} &= rC_1 \,,
\\
M_{(2)} &= r(r+1)C_1^2 + rC_2 \,,
\\
\sigma^2 &= M_{(2)} - M^2_{(1)} +M_{(1)}
\nonumber\\
&= r(r+1)C_1^2 + rC_2 -r^2C_1^2 +rC_1
\nonumber\\
&= r(C_2 +C_1^2 + C_1) \,.
\end{align}
\end{subequations}
\item
  Similarly, using $F$ indexing,
\begin{subequations}
\begin{align}
\hat{M}_{(1)} &= rF_1 \,,
\\
\hat{M}_{(2)} &= r(r+1)F_1^2 + rF_2 \,,
\\
\sigma^2 &= \hat{M}_{(2)} - \hat{M}^2_{(1)} +\hat{M}_{(1)}
\nonumber\\
&= r(r+1)F_1^2 + rF_2 -r^2F_1^2 +rF_1
\nonumber\\
&= r(F_2 +F_1^2 + F_1) \,.
\end{align}
\end{subequations}
\item
  The expression for the variance is well known and is
\bq  
\sigma^2 = r\,\biggl[\frac{1}{(qp^k)^2} -\frac{2k+1}{qp^k} -\frac{p}{q^2}\biggr] \,.
\eq
\end{enumerate}

\newpage
\subsection{Summary of results}
Let us collect the results in one place, for ease of reference and also to establish the similiarities of the expressions using the two indexing schemes.
The factorial means are given by the same formal sum, just substitute $C_j$ for $C$ indexing and $F_j$ for $F$ indexing.
Recall $s_n = n_1+\dots+n_k$.
\begin{subequations}
\label{eq:facmom_summary}
\begin{align}
M_{(n)} = n! \sum_{n_1+2n_2+\dots+kn_k=n} (r+s_n-1)_{s_n} \prod_{j=1}^k \frac{1}{n_j!}\Bigl(\frac{C_j}{j!}\Bigr)^{n_j} \,,
\\
\hat{M}_{(n)} = n! \sum_{n_1+2n_2+\dots+kn_k=n} (r+s_n-1)_{s_n} \prod_{j=1}^k \frac{1}{n_j!}\Bigl(\frac{F_j}{j!}\Bigr)^{n_j} \,.
\end{align}
\end{subequations}
The expressions for $C_j$ and $F_j$ are 
\begin{subequations}
\label{eq:Cj_Fj_summary}
\begin{align}
  C_j &= \frac{j!}{q^jp^k}\biggl[ p^{j-1}(1-p^k) -p^k\sum_{i=1}^{j-1} \binom{k}{i} p^{j-i-1}q^i \biggr] -(k)_j \,,
  \\
  F_j &= \frac{j!}{q^jp^k}\biggl[ p^{j-1}(1-p^k) +\sum_{i=1}^{j-1} (-1)^i \binom{k+i-1}{i} p^{j-i-1}q^i \biggr] \,.
\end{align}
\end{subequations}
Some of the terms in $C_j$ and $F_j$ can be collected into hypergeometric functions
\begin{subequations}
\label{eq:Cj_Fj_hyp_summary}
\begin{align}
  C_j &= \frac{j!}{q^jp^k}\biggl[ p^{j-1} -\binom{k}{j-1}p^kq^{j-1}\,{}_2F_1(1-j,1,k-j+2;-p/q) \biggr] -(k)_j \,,
  \\
  F_j &= \frac{j!}{q^jp^k}\biggl[ -p^{k+j-1} +(-1)^{j-1}\binom{k+j-2}{j-1}q^{j-1}\,{}_2F_1(1-j,1,2-j-k;-p/q) \biggr] \,.
\end{align}
\end{subequations}

\newpage
\setcounter{equation}{0}
\section{\label{sec:mom}Raw (power) moments}
\subsection{General remarks}
Denote the raw (or power) moment by $M_n$ for $C$ indexing and by $\hat{M}_n$ for $F$ indexing.
The subscript is $n$ not $(n)$.
Again, let $Y \sim \mathrm{NB}_{k,r}$ with $F$ indexing and $Z \sim \mathrm{NB}_{k,r}$ with $C$ indexing.
Then by definition
$\hat{M}_n = \mathbb{E}[Y^n]$ and $M_n = \mathbb{E}[Z^n]$ and it is also evident that
\bq
Y^n = (Z+rk)^n = \sum_{i=0}^n \binom{n}{i}Z^{n-i} (rk)^i \,.
\eq
Hence
\bq
\label{es:mom_shift}
\hat{M}_n = \sum_{i=0}^n \binom{n}{i}M_{n-i} (rk)^i \,.
\eq
We follow the ideas in \cite{Mane_CC23_12}.
First define functions $\tilde{C}_j$ for $C$ indexing and $\tilde{F}_j$ for $F$ indexing.
The expressions for the raw moments $M_n$ and $\hat{M}_n$ are formally similar to those for the factorial moments in eq.~\eqref{eq:facmom_summary}.
We substitute $\tilde{C}_j$ in place of $C_j$ or $\tilde{F}_j$ in place of $F_j$ to obtain the sum for the raw moments.
Recall $s_n = n_1+\dots+n_k$.
\begin{subequations}
\label{eq:facmom_summary}
\begin{align}
M_n = n! \sum_{n_1+2n_2+\dots+kn_k=n} (r+s_n-1)_{s_n} \prod_{j=1}^k \frac{1}{n_j!}\Bigl(\frac{\tilde{C}_j}{j!}\Bigr)^{n_j} \,,
\\
\hat{M}_n = n! \sum_{n_1+2n_2+\dots+kn_k=n} (r+s_n-1)_{s_n} \prod_{j=1}^k \frac{1}{n_j!}\Bigl(\frac{\tilde{F}_j}{j!}\Bigr)^{n_j} \,.
\end{align}
\end{subequations}
It remains to write explicit expressions for $\tilde{C}_j$ and $\tilde{F}_j$.
They have a pattern, but it is more complicated than for the factorial moments.
They involve the Eulerian polynomials, as will be explained below.

The Eulerian polynomials $A_i(t)$ are defined as the coefficients of $x^i$ as follows
\bq
\sum_{i=0}^\infty A_i(t)\frac{x^i}{i!} = \frac{t-1}{t-e^{(t-1)x}} \,.
\eq
Here $A_0(t)=1$ and for $i\ge1$, $A_i(t)$ is a polynomial in $t$ of degree $i-1$.
Also $A_i(t) = \sum_{j=0}^{i-1}A_{ij}t^j$ and
$A_{ij}$ is an element in the triangle of Eulerian Numbers, given by
\bq
A_{ij} = \sum_{s=0}^j(-1)^s\binom{i+1}{s}(j+1-s)^i \,. 
\eq
The first few Eulerian polynomials are  
\begin{subequations}
\label{eq:Ei_poly_ex}
\begin{align}
  A_0(t) &= 1 \,,
  \\
  A_1(t) &= 1 \,,
  \\
  A_2(t) &= 1+t \,,
  \\
  A_3(t) &= 1+4t+t^2 \,,
  \\
  A_4(t) &= 1+11t+11t^2+t^3 \,.
\end{align}
\end{subequations}

\newpage
\subsection{$C$ indexing}
We begin with $C$ indexing.
Define $\xi_j$ via
\bq
\tilde{C}_j = \frac{\xi_j}{q^jp^k} \,.
\eq
Then we obtain the following expressions for $\xi_j$ for $j=1,\dots,4$.
\begin{subequations}
\label{eq:xi_j_ex}
\begin{align}
\xi_1 &= (1-p^k)A_1(p) -kp^k \,,
\\
\xi_2 &= (1-p^k)A_2(p) -p^k\Bigl[ 2kqA_1(p) +(kq)^2\Bigr] \,,
\\
\xi_3 &= (1-p^k)A_3(p) -p^k\Bigl[ 3kqA_2(p) +3(kq)^2A_1(p) +(kq)^3\Bigr] \,,
\\
\xi_4 &= (1-p^k)A_4(p) -p^k\Bigl[ 4kqA_3(p) +6(kq)^2A_2(p) +4(kq)^3A_1(p) +(kq)^4 \Bigr] \,.
\end{align}
\end{subequations}
In general, the expression is
\bq
\label{eq:xi_j}
\xi_j = (1-p^k)A_j(p) -p^k\sum_{i=1}^j \binom{j}{i}(kq)^iA_{j-i}(p) \,.
\eq

\newpage
\subsection{$F$ indexing}
The pattern for $F$ indexing is similar, but the signs alternate and the last term
(highest power of $kq$) is absent and there is no factor of $p^k$, as will be explained below.
Define $\chi_j$ via
\bq
\tilde{F}_j = \frac{\chi_j}{q^jp^k} \,.
\eq
Then we obtain the following expressions for $\chi_j$ for $j=1,\dots,4$.
\begin{subequations}
\label{eq:chi_j_ex}
\begin{align}
\chi_1 &= (1-p^k)A_1(p) \,,
\\
\chi_2 &= (1-p^k)A_2(p) -2kqA_1(p) \,,
\\
\chi_3 &= (1-p^k)A_3(p) -3kqA_2(p) +3(kq)^2A_1(p) \,,
\\
\chi_4 &= (1-p^k)A_4(p) -4kqA_3(p) +6(kq)^2A_2(p) -4(kq)^3A_1(p) \,.
\end{align}
\end{subequations}
In general, the expression is
\bq
\label{eq:chi_j}
\chi_j = (1-p^k)A_j(p) +\sum_{i=1}^{j-1} (-1)^i\binom{j}{i}(kq)^iA_{j-i}(p) \,.
\eq

\newpage
\setcounter{equation}{0}
\section{\label{sec:cent}Central moments}
\subsection{General remarks}
The term ``central moments'' refers to moments centered on the mean.
The central moments do not depend on the indexing scheme.
Let $Y \sim \mathrm{NB}_{k,r}$ with $F$ indexing and $Z \sim \mathrm{NB}_{k,r}$ with $C$ indexing.
Denote their respective means by $\mu_F$ and $\mu_C$ and recall $\mu_F = \mu_C +rk$.
Then
\bq
\mathbb{E}[(Y-\mu_F)^n] = \mathbb{E}[(Z+rk-(\mu_C+rk))^n] = \mathbb{E}[(Z-\mu_C)^n] \,.
\eq
Denote the central moments by $\mathcal{M}_n$. By definition $\mathcal{M}_1=0$ and $\mathcal{M}_2$ is the variance $\sigma^2$.
\bq  
\begin{split}
  \sigma^2 = \mathcal{M}_2 &= M_2 -\mu_C^2
  \\
  &= r(r+1)\tilde{C}_1^2 +r\tilde{C}_2 -r^2\tilde{C}_1^2
  \\
  &= r(\tilde{C}_1^2 +\tilde{C}_2)
  \\
  &= r(\tilde{F}_1^2 +\tilde{F}_2) \,.
\end{split}
\eq
As already noted, the variance is proportional to $r$.

\newpage
\subsection{Third central moment}
The third central moment is given as follows.
\bq  
\begin{split}
  \mathcal{M}_3 &= \mathbb{E}[(Z - \mu_C)^3]
  \\
  &= M_3 -3\mu_C M_2 +3\mu_C^2M_1 -\mu_C^3
  \\
  &= r(r+1)(r+2)\tilde{C}_1^3 +3r(r+1)\tilde{C}_1\tilde{C}_2 +r\tilde{C}_3 -3r\tilde{C}_1\Bigl[r(r+1)\tilde{C}_1^2+r\tilde{C}_2\Bigr] +3r^3\tilde{C}_1^3 -r^3\tilde{C}_1^3
  \\
  &= r(2\tilde{C}_1^3 +3\tilde{C}_1\tilde{C}_2 +\tilde{C}_3)
  \\
  &= r(2\tilde{F}_1^3 +3\tilde{F}_1\tilde{F}_2 +\tilde{F}_3) \,.
\end{split}
\eq
The third central moment is proportional to $r$.
The skewness $\gamma$ is the normalized third central moment 
\bq
\gamma = \frac{\mathcal{M}_3}{\sigma^3} \propto r^{-1/2} \,.
\eq
For fixed $k$ and $p$, the skewness {\em decreases} as $r$ increases.

\newpage
\subsection{Fourth central moment}
The fourth central moment is given as follows.
\bq  
\begin{split}
  \mathcal{M}_4 &= \mathbb{E}[(Z - \mu_C)^4]
  \\
  &= M_4 -4\mu_C M_3 +6\mu_C^2M_2 -4\mu_C^3M_1 +\mu_C^4
  \\
  &= r(r+1)(r+2)(r+3)\tilde{C}_1^4 +6r(r+1)(r+2)\tilde{C}_1^2\tilde{C}_2 +r(r+1)(4\tilde{C}_1\tilde{C}_3+3\tilde{C}_2^2) +r\tilde{C}_4
  \\
  &\qquad
  -4r\tilde{C}_1\Bigl[r(r+1)(r+2)\tilde{C}_1^3 +3r(r+1)\tilde{C}_1\tilde{C}_2 +r\tilde{C}_3\Bigr]
  \\
  &\qquad  
  +6r^2\tilde{C}_1^2\Bigl[r(r+1)\tilde{C}_1^2+r\tilde{C}_2\Bigr] -4r^4\tilde{C}_1^4 +r^4\tilde{C}_1^4
  \\
  &= 3r^2(\tilde{C}^4 +2\tilde{C}_1^2\tilde{C}_2 +\tilde{C}_2^2)
  +r(6\tilde{C}^4 +12\tilde{C}_1^2\tilde{C}_2 +3\tilde{C}_2^2 +4\tilde{C}_1\tilde{C}_3 +\tilde{C}_4)
  \\
  &= 3\sigma^4
  +r(6\tilde{C}^4 +12\tilde{C}_1^2\tilde{C}_2 +3\tilde{C}_2^2 +4\tilde{C}_1\tilde{C}_3 +\tilde{C}_4) 
  \\
  &= 3\sigma^4
  +r(6\tilde{F}^4 +12\tilde{F}_1^2\tilde{F}_2 +3\tilde{F}_2^2 +4\tilde{F}_1\tilde{F}_3 +\tilde{F}_4) \,.
\end{split}
\eq
The fourth central moment has terms of both $O(r^2)$ and $O(r)$.
However, for the {\em excess kurtosis} we subtract $3\sigma^4$, which leaves only terms proportional to $r$.
\bq
\frac{\mathcal{M}_4-3\sigma^4}{\sigma^4} \propto r^{-1} \,.
\eq
For fixed $k$ and $p$, the excess kurtosis {\em decreases} as $r$ increases.

\newpage
\setcounter{equation}{0}
\section{\label{sec:DM4}Alternative (non-combinatorial) viewpoint}
We employ $F$ indexing exclusively in this section.
For a sequence of independent Bernoulli trials with probability $p\in(0,1)$ of success,  
let $N_k$ be the waiting time for the first run of $k$ consecutive successes. 
Then $N_k$ is said to have the \emph{geometric distribution of order $k$}.  
For $k=1$ it is the geometric distribution.
Its probability generating function (pgf) is \cite{Feller} 
\bq
\label{eq:pgf}
\phi(s,k) = \frac{p^ks^k(1-ps)}{1-s+qp^ks^{k+1}} = \frac{p^ks^k}{1-qs(1+pk+\dots+k^{k-1}s^{k-1})} \,.
\eq
The pmf $f_k(n)$ satisfies the recurrence \cite{BarryLoBello}
\bq
\label{eq:recrel}
f_k(n) = qf_k(n-1) + pqf_k(n-2) + p^2qf_k(n-3) +\dots + p^{k-1}qf_k(n-k) \,.
\eq
The initial conditions are $f_k(n)=0$ for $n\in[1,k-1]$ and $f_k(k)=p^k$.
The auxiliary polynomial is
\bq
\label{eq:auxpoly}
\mathscr{A}(z) = z^k -qz^{k-1} -qpz^{k-2} -\dots -qp^{k-1} \,.
\eq
Here $z=1/s$ (see eq.~\eqref{eq:pgf}).
Feller \cite{Feller} proved the roots of the auxiliary polynomial are distinct and have magnitude less than unity.
We denote the roots by $\lambda_j(p,k)$, where $j=0,\dots,k-1$.
An important fact is that for any root, $\lambda_j^k(1-\lambda_j) = p^kq$ (see eq.~\eqref{eq:auxpoly} or \cite{DM3}).
Dilworth and Mane \cite{DM3} derived expressions for the roots in terms of Fuss-Catalan numbers.
The details can be found in \cite{DM3} and will not be repeated here.

In a second paper, Dilworth and Mane \cite{DM4} applied the formalism in \cite{DM3}
to treat the negative binomial distribution NB$(k,r)$ of Type I with $\ell$-overlapping, say NB$(k,r,\ell)$.
Up to now, we have considered only the case $\ell=0$ of nonoverlapping success runs.
The results below apply for all $\ell \in [0,k-1]$.
The pgf for NB$(k,r,\ell)$ is (eq.~(24) in \cite{DM4})
\bq
\phi_{r,k,\ell}(s) = \frac{\phi(s,k)^r}{\phi(s,\ell)^{r-1}} \,.
\eq
By definition, the pgf $\phi_{r,k,\ell}(s)$ is related to the pmf $f_{r,k,\ell}(n)$ via
\bq
\phi_{r,k,\ell}(s) = \sum_{n=0}^\infty s^n f_{r,k,\ell}(n) \,.
\eq
Then the pmf is given as a sum of the powers of the roots $\lambda_j$ as follows (eq.~(52) in \cite{DM4})
\bq
f_{r,k,\ell}(n) = \sum_{j=0}^{k-1}\sum_{m=1}^r\binom{n-1}{m-1}a_{jm}\lambda_j(p,k)^{n-m} \,.
\eq
The details of the coefficients $a_{jm}$ were given in \cite{DM4} and their explicit values are not required below.

To avoid confusion with too many uses of $n$, we denote the $i^{th}$ factorial moment by $\mathscr{M}_{(i)}$.
Then
\bq
\label{eq:facmom_DM4}
\begin{split}
  \mathscr{M}_{(i)} &= i! \sum_{j=0}^{k-1}\sum_{m=1}^r a_{jm}\lambda_j^{-m} \biggl[ \sum_{n=i}^\infty \binom{n}{i} \binom{n-1}{m-1}\lambda_j^n \biggr]
  \\
  &= i! \sum_{j=0}^{k-1}\sum_{m=1}^r a_{jm}\lambda_j^{-m} \biggl[\frac{\lambda_j^i}{(1-\lambda_j)^{i+m}}\,\sum_{b=1}^m \binom{m}{b}\binom{i-1}{b-1}\lambda_j^{m-b}\biggr]
  \qquad \textrm{(because $|\lambda_j| < 1$)}
  \\
  &= i! \sum_{j=0}^{k-1}\sum_{m=1}^r a_{jm} \frac{\lambda_j^i}{(1-\lambda_j)^{i+m}}\,\sum_{b=1}^m \binom{m}{b}\binom{i-1}{b-1}\lambda_j^{-b}
  \\
  &= i! \sum_{j=0}^{k-1}\sum_{m=1}^r a_{jm} \frac{\lambda_j^{(i+m)k+i}}{\lambda_j^{(i+m)k}(1-\lambda_j)^{i+m}}\,\sum_{b=1}^m \binom{m}{b}\binom{i-1}{b-1}\lambda_j^{-b}
  \\
  &= i! \sum_{j=0}^{k-1}\sum_{m=1}^r a_{jm} \frac{\lambda_j^{(i+m)k+i}}{(qp^k)^{i+m}}\,\sum_{b=1}^m \binom{m}{b}\binom{i-1}{b-1}\lambda_j^{-b}
  \qquad\qquad\qquad \textrm{(because $\lambda_j^k(1-\lambda_j) = p^kq$)}
  \\
  &= \frac{i!}{(qp^k)^i} \sum_{j=0}^{k-1}\sum_{m=1}^r a_{jm} \frac{\lambda_j^{(i+m)k+i}}{(qp^k)^m}\,\sum_{b=1}^m \binom{m}{b}\binom{i-1}{b-1}\lambda_j^{-b}
  \\
  &= \frac{i!}{(qp^k)^i} \sum_{j=0}^{k-1}\sum_{m=1}^r a_{jm} \frac{m\lambda_j^{(i+m)k+i-1}}{(qp^k)^m}\,{}_2F_1(1-m,1-i,2;1/\lambda_j) \,.
\end{split}
\eq
In the special case $r=1$, then $m=1$ hence $b=1$ (or ${}_2F_1(0,\dots)=1$) and the above simplifies to the expression derived in \cite{Mane_CC23_13}
\bq
\begin{split}
  \mathscr{M}_{(i)} &= \frac{i!}{(qp^k)^{i+1}} \sum_{j=0}^{k-1} a_{j1} \lambda_j^{(i+1)k+i} \lambda_j^{-1}
  \\
  &= \frac{i!}{(qp^k)^{i+1}} \,f_k((i+1)k+i) \,.
\end{split}
\eq
Note the following.
\begin{enumerate}
\item
  First, eq.~\eqref{eq:facmom_DM4} furnishes a non-combinatorial expression for the factorial moments.
\item
  Moreover, it does so not only for nonoverlapping success runs but for {\em all} $\ell \in [0,k-1]$.
\item
  Note that the sum in eq.~\eqref{eq:facmom_DM4} is {\em independent of $n$}.
  It contains $rk$ summations (and the computation of the hypergeometric function).
\item
  Of course, one must calculate the roots $\lambda_j$ and the coefficients $a_{jm}$, which may require a numerical algorithm.
\end{enumerate}
Finally, note that the case $\ell=k-1$ is the negative binomial distribution of Type III.
Hence eq.~\eqref{eq:facmom_DM4} yields the factorial moments for the negative binomial distribution of Type III.

\newpage
\setcounter{equation}{0}
\section{\label{sec:gap}Gaps between success runs}
Most of the notation and formalism below is based on the exposition in Sec.~\ref{sec:DM4}.
It was noted in \cite{DM4} that the analysis therein also applies to the case where there is a minimum gap between the success runs.
This can be modelled as a ``negative overlap'' which means $\ell<0$ (see the presentation in Sec.~\ref{sec:DM4} above).
We must clarify what ``gap'' means in this context.
Denote the gap by $g$, where $g>0$.
Then after a success run of length $k$ is attained, {\em the outcomes of the subsequent $g$ trials are ignored}.
The count for the next success runs begins only after those $g$ trials have been completed.
To return to the example in Sec.~\ref{sec:intro}, suppose $k=2$ and $g=2$.
Then the sequence $SSSSSS$ contains two success runs of length 2, because the outcomes of the third and fourth trials are ignored.
Similarly, the sequences $SSFFSS$, $SSSFSS$ and $SSFSSS$ all contain two success runs of length 2.
In all four cases, the two success runs are the first two and last two trials.
This is {\em not} the same as the negative binomial distribution of Type II.
The latter requires at least one failure between success runs. The number of trials between success runs is arbitrary (but at least one).

We employ $F$ indexing exclusively in this section.
With respect to the formalism in \cite{DM4} and the exposition in Sec.~\ref{sec:DM4}, $g=-\ell$ (and $g>0$).
Then it was noted in \cite{DM4} that the pmf for $\ell<0$ is given by that for $\ell=0$ via the relation (Prop.~19 in \cite{DM4})
\bq
\label{eq:pmg_gap}
f_{r,k,\ell}(n) = f_{r,k,0}(n +(r-1)\ell) = f_{r,k,0}(n -(r-1)g) \,.
\eq
Hence the pgf is given by
\bq
\phi_{r,k,-g}(s) = \phi_{r,k,0}(s)s^{(r-1)g} \,.
\eq
The mean and variance are
\begin{subequations}
\begin{align}
\mu_{r,k,-g} &= \mu_{r,k,0} +(r-1)g &=&\; r\mu_{1,k,0} +(r-1)g \,,
\\
\sigma^2_{r,k,-g} &= \sigma^2_{r,k,0} &=&\; r\sigma^2_{1,k,0} \,.
\end{align}
\end{subequations}
The expressions for the pgf and mean correct misprints in \cite{DM4}.
Next, eq.~\eqref{eq:pmg_gap} furnishes the key to computing the moments and factorial moments for $G$.
We reuse the notation $\mathscr{M}_{(i)}$ from Sec.~\ref{sec:DM4} but attach arguments $(r,k,\ell) = (r,k,-g)$ to clarify the dependency on its arguments.
We obtain, analogous to eqs.~\eqref{es:facmom_shift} and eq.~\eqref{es:mom_shift} for the factorial moments and raw moments respectively,
\begin{subequations}
\begin{align}
\mathscr{M}_{(n)}(r,k,-g) &= \sum_{i=0}^n \binom{n}{i} \hat{M}_{(n-i)}((r-1)g)_{(i)} \,,
\\
\mathscr{M}_{n}(r,k,-g) &= \sum_{i=0}^n \binom{n}{i} \hat{M}_{n-i}((r-1)g)^i \,.
\end{align}
\end{subequations}
The central moments are the same as those for a gap of zero, as was noted above for the variance.

\newpage
\setcounter{equation}{0}
\section{\label{sec:conc}Conclusion}
Most of this note was devoted to the negative binomial distribution NB$(k,r)$ of Type I,
which is the probability distribution for a sequence of independent Bernoulli trials
(with success parameter $p\in(0,1)$) with $r$ nonoverlapping success runs of length $\ge k$.
We presented a new, more concise, expression for its probability mass function,
and also noted that our expression can be succinctly written using hypergeometric functions.
In this context, the published literature employs two different indexing schemes to count the trials,
hence care was taken to distinguish between them and present expressions for both schemes.
We also presented new expressions (as combinatorial sums) for the moments and factorial moments, as opposed to only the mean and variance (which are already known).
We also briefly discussed the central moments, displaying in particular expressions for the skewness and kurtosis.
Next, we presented an alternative non-combinatorial viewpoint, where the probability mass function was expressed as
a sum of powers of the roots of the auxiliary equation associated to the distribution.
This non-combinatorial formalism allowed us to derive expressions for the factorial moments not only for nonoverlapping success runs,
but also for runs with an overlap of $\ell$, where $\ell\in[0,k-1]$.
The case $\ell=k-1$ is the negative binomial distribution NB$(k,r)$ of Type III.
The non-combinatorial formalism also treats the scenario of a minimum gap between success runs, effectively a `negative overlap' $g=-\ell > 0$.
Expressions for the moments and factorial moments were presented.

\newpage
\setcounter{equation}{0}
\section{\color{red}Addendum 1/23/2024}
We can derive the probability mass function directly from the probability generating function in eq.~\eqref{eq:pgf}.
We employ $F$ indexing. 
We use the binomial theorem to expand the pgf in powers of $s$ and temporarily set the upper limits on the sums to $\infty$ and determine their exact values later.
The pgf in this case is
\bq
\begin{split}
  \psi(s,k) = \phi(s,k)^r &= \frac{p^{rk}s^{rk}(1-ps)^r}{(1-s)^r}\frac{1}{(1 +qp^ks^{k+1}/(1-s))^r}
  \\
  &= p^{rk}s^{rk}(1-ps)^r \sum_{i=0}^\infty (-1)^i\binom{r+i-1}{r-1}\frac{q^ip^{ik}s^{i(k+1)}}{(1-s)^{r+i}}
  \\
  &= p^{rk}s^{rk}(1-ps)^r \sum_{i=0}^\infty (-1)^i\binom{r+i-1}{r-1}q^ip^{ik}s^{i(k+1)}\sum_{u=0}^\infty \binom{r+i+u-1}{r+i-1} s^u
  \\
  &= \dots
  \\
  &= p^{rk}s^{rk}\sum_{v=0}^\infty s^v \sum_{j=0}^\infty (-1)^{v-j}p^{v-j} \binom{r}{v-j} \sum_{i=0}^\infty (-1)^iq^ip^{ik} \binom{r+i-1}{r-1}\binom{r+j-ik-1}{r+i-1} \,.
\end{split}
\eq
The probability $P_n$ is the coefficient of $s^n$, hence set $n = v+rk$ or $v = n-rk$.
\bq
  P_n = p^{rk} \sum_{j=0}^{j_{\max}}(-1)^{n-rk-j}p^{n-rk-j} \binom{r}{n-rk-j} \sum_{i=0}^{i_{\max}} (-1)^iq^ip^{ik} \binom{r+i-1}{r-1}\binom{r+j-ik-1}{r+i-1} \,.
\eq
The value of $i_{\max}$ is given by $j-ik=i$ hence $i \le \lfloor j/(k+1) \rfloor$.  
The value of $j_{\max}$ is given by $0 = n - rk - j$ hence $j \le n - rk$.
Then
\bq
\label{eq:Pn_from_pgf}
  P_n = p^{rk} \sum_{j=0}^{n-rk}(-1)^{n-rk-j}p^{n-rk-j} \binom{r}{n-rk-j} \sum_{i=0}^{\lfloor j/(k+1) \rfloor} (-1)^iq^ip^{ik} \binom{r+i-1}{r-1}\binom{r+j-ik}{r+i} \,.
\eq
The sums in eq.~\eqref{eq:Pn_from_pgf} are a double nested sum,
as opposed to the multinomial sum derived by 
Philippou, Georghiou, and Philippou (Theorem 3.1 in \cite{PhilippouGeorghiouPhilippou}),
who employed the multinomial theorem to expand the pgf in eq.~\eqref{eq:pgf}.
The sums in eq.~\eqref{eq:Pn_from_pgf} are not as concise but are equivalent to eq.~\eqref{eq:pmf_sum_F}.
The factorial moment $\hat{M}_{(n)}$ is given by the $n^{th}$ partial derivative of the pgf via
\bq
\label{eq:facmom_from_pgf}
\begin{split}
\hat{M}_{(n)} &= \biggl[\frac{\partial^n\psi(s,k)}{\partial s^n}\biggr]_{s=1}
\\
&= n!p^{rk}\sum_{v=n-rk}^\infty \binom{v+rk}{n} \sum_{j=0}^{v-rk} (-1)^{v-j}p^{v-j} \binom{r}{v-j} \sum_{i=0}^{\lfloor j/(k+1) \rfloor} (-1)^iq^ip^{ik} \binom{r+i-1}{r-1}\binom{r+j-ik-1}{r+i-1} \,.
\end{split}
\eq
This does not look as elegant as eq.~\eqref{eq:facmom_F}, but is a triple nested sum,
as opposed to the more complicated $k$-fold nesting in eq.~\eqref{eq:facmom_F}.
However, the sum over $v$ in eq.~\eqref{eq:facmom_from_pgf} is an infinite sum,
whereas the sums and products in eq.~\eqref{eq:facmom_F} are finite.

\newpage
\setcounter{equation}{0}
\section{\color{red}Addendum 1/26/2024}
This section comments on Muselli's paper \cite{Muselli} on the negative binomial distribution NB$(k,r)$ of Type II.
As noted in eq.~\eqref{eq:Muselli_pmf}, Muselli published elegant expressions for success runs in Bernoulli trials.
We comment on two theorems in \cite{Muselli}.
They are expressed below in terms of the notation employed in this note.
\begin{enumerate}
\item
  Muselli Theorem 3 states that if $M_n^{(k)}$ is the number of success runs with length $k$ or more in $n$ Bernoulli trials, then
  (where $k$, $n$ and $r$ are positive integers)
\bq
\label{eq:Muselli_Thm3}
P(M_n^{(k)}=r) = \sum_{m=r}^{\left\lfloor\frac{n+1}{k+1}\right\rfloor} (-1)^{m-r}\binom{m}{r} p^{mk}q^{m-1} \biggl[\binom{n-mk}{m-1} +q\binom{n-mk}{m}\biggr] \,.
\eq
\item
  Muselli Theorem 4 states that the random variable NB${}_{k,r}$ is characterized by the following probability distribution (see eq.~\eqref{eq:Muselli_pmf} above)
\bq
\label{eq:Muselli_Thm4}
P(\textrm{NB}_{k,r}=n) = \sum_{m=r}^{\lfloor\frac{n+1}{k+1}\rfloor} (-1)^{m-r}\binom{m-1}{r-1}p^{mk}q^{m-1}\biggl[\binom{n-mk-1}{m-2} +q\binom{n-mk-1}{m-1}\biggr] \,.
\eq
\end{enumerate}
It was noted in \cite{Mane_CC23_13} that Muselli's expressions contain difficulties of interpretation.
Muselli stated that $r\ge1$ in eq.~\eqref{eq:Muselli_Thm3}, but in fact we must allow $r=0$.
(This is required to derive Muselli's Theorem 4, because Theorem 4 requires eq.~\eqref{eq:Muselli_Thm3} with $r-1$, and $r-1=0$ if $r=1$.)
Setting $r=0$ in eq.~\eqref{eq:Muselli_Thm3}, then the sum
in eq.~\eqref{eq:Muselli_Thm3} begins at $m=0$ and we obtain $\binom{n-mk}{m-1} = \binom{n}{-1}$.
In eq.~\eqref{eq:Muselli_Thm4}, we must have $r\ge1$. 
If $r=1$, then the sum in eq.~\eqref{eq:Muselli_Thm4} begins at $m=1$ and we obtain $\binom{n-mk-1}{m-2} = \binom{n-k-1}{-1}$.
We rewrite Muselli's expressions as follows, employing the identity
\bq
\binom{i-1}{j-1} = \binom{i}{j} - \binom{i-1}{j} \,.
\eq
\begin{enumerate}
\item
  We reexpress eq.~\eqref{eq:Muselli_Thm3} as follows.
\bq
\label{eq:Muselli_Thm3_alt}
\begin{split}
  P(M_n^{(k)}=r) &= \sum_{m=r}^{\left\lfloor\frac{n+1}{k+1}\right\rfloor} (-1)^{m-r}\binom{m}{r} p^{mk}q^{m-1} 
  \biggl[\binom{n-mk+1}{m} -\binom{n-mk}{m} +q\binom{n-mk}{m}\biggr]
  \\
  &= \sum_{m=r}^{\left\lfloor\frac{n+1}{k+1}\right\rfloor} (-1)^{m-r}\binom{m}{r} p^{mk}q^{m-1} \biggl[\binom{n-mk+1}{m} -p\binom{n-mk}{m}\biggr] \,.
\end{split}
\eq
\item
  We reexpress eq.~\eqref{eq:Muselli_Thm4} as follows.
\bq
\label{eq:Muselli_Thm4_alt}
\begin{split}
  P(\textrm{NB}_{k,r}=n) &= \sum_{m=r}^{\lfloor\frac{n+1}{k+1}\rfloor} (-1)^{m-r}\binom{m-1}{r-1}p^{mk}q^{m-1} 
  \biggl[\binom{n-mk}{m-1} -\binom{n-mk-1}{m-1} +q\binom{n-mk-1}{m-1}\biggr]
  \\
  &= \sum_{m=r}^{\lfloor\frac{n+1}{k+1}\rfloor} (-1)^{m-r}\binom{m-1}{r-1}p^{mk}q^{m-1} \biggl[\binom{n-mk}{m-1} -p\binom{n-mk-1}{m-1}\biggr] \,.
\end{split}
\eq
\end{enumerate}
Then eqs.~\eqref{eq:Muselli_Thm3_alt} and \eqref{eq:Muselli_Thm4_alt} do not contain binomial coefficients with negative denominators.

\end{document}